\documentclass[conference]{IEEEtran}
\IEEEoverridecommandlockouts
\usepackage{cite}
\usepackage{amsmath,amssymb,amsfonts}
\usepackage{graphicx}
\usepackage{textcomp}
\usepackage{xcolor}
\def\BibTeX{{\rm B\kern-.05em{\sc i\kern-.025em b}\kern-.08em
    T\kern-.1667em\lower.7ex\hbox{E}\kern-.125emX}}

\usepackage{mathrsfs}
\usepackage{multicol}
\usepackage{wrapfig}
\usepackage[utf8]{inputenc}
\usepackage{amsthm}
\usepackage[shortlabels]{enumitem}
\usepackage{here}
\usepackage{multirow}
\usepackage{stmaryrd}
\usepackage{version}
\usepackage{comment}
\usepackage{listings}
\usepackage{verbatim}
\usepackage{xcolor}
\usepackage{enumitem}
\usepackage{cases}
\usepackage{xspace}
\usepackage{algorithm2e}
\RestyleAlgo{ruled}
\SetAlCapNameFnt{\small}
\SetAlCapFnt{\small}
\usepackage{algpseudocode}
\newcommand{\1}{\mbox{\fontencoding{U}\fontfamily{bbold}\selectfont1}}

\DeclareMathAlphabet{\mathcal}{OMS}{cmsy}{m}{n}
\newtheorem{theorem}{Theorem}

\newtheorem{prop}[theorem]{Proposition}
\newtheorem{defn}{Definition}

\usepackage{hyperref}%
\hypersetup{colorlinks,allcolors=black}

\newcommand{\E}{\mathbb{E}}
\newcommand{\given}{\,\vert\,}
\DeclareMathOperator*{\argmin}{arg\,min}
\DeclareMathOperator*{\argmax}{arg\,max}

\newcommand{\Hsquare}{%
  \text{\fboxsep=-.2pt\fbox{\rule{0pt}{1ex}\rule{1ex}{0pt}}}%
}

\newcounter{hints}
\renewcommand{\thehints}{\roman{hints}}
\newcommand{\hintedrel}[2][]{%
  \stepcounter{hints}%
  \if\relax\detokenize{#1}\relax\else\csxdef{hint@#1}{\thehints}\fi
  \mathrel{\overset{\textrm{(\thehints)}}{\vphantom{\le}{#2}}}%
}
\newcommand{\hintref}[1]{\csuse{hint@#1}}

\begin{document}

\title{
The Price of Simplicity: Analyzing Decoupled Policies for Multi-Location Inventory Control
% The Price of Simplicity: Performance Guarantees for Decoupled Policies in Multi-Location Inventory Control
%Analysis of Decoupled Policies for Multi-Location Inventory Control Problems

\author{
Yohan John*, Vade Shah*, James A. Preiss, Mahnoosh Alizadeh, and Jason R. Marden

\thanks{This paper was supported by the California Energy Commission, project \#EPIC-24-012.}
\thanks{Y. John ({\tt\small yohanjohn@ucsb.edu}), V. Shah, J. A. Preiss, M. Alizadeh, and J. R. Marden are with the Center for Control, Dynamical Systems, and Computation at the University of California, Santa Barbara, CA. * denotes equal contribution.}%
}

}

\maketitle

\begin{comment}
\begin{abstract}
    In multidimensional control systems, optimal control policies may be complex functions of the entire state space. Single-dimensional decoupled policies offer a simpler alternative, but their simplicity may come at the expense of optimality. In this work, we study the performance gap between simple decoupled policies and optimal coupled policies in the context of the well-known inventory control problem, where a planner seeks to identify optimal inventory levels that meet stochastic demand while minimizing costs associated with ordering additional inventory and carrying excess or insufficient inventory. Inspired by ideas of `economies of scale', we consider multi-location systems where the ordering cost is a nonlinear function of the cumulative order amount across all locations. Our main results characterize the worst-case performance guarantees of simple decoupled policies in broad classes of coupled systems, and a key takeaway is that these guarantees depend highly upon the nonlinearity of the cost function. Numerical simulations demonstrate that although decoupled policies may be suboptimal, they significantly outperform their worst-case bounds for typical problem instances. These findings provide practical guidance for navigating the trade-off between control complexity and operational efficiency in a wide range of coupled systems. 
\end{abstract}
\end{comment}

\begin{abstract}
    What is the performance cost of using simple, decoupled control policies in inherently coupled systems? Motivated by industrial refrigeration systems, where centralized compressors exhibit economies of scale yet traditional control employs decoupled room-by-room temperature regulation, we address this question through the lens of multi-location inventory control. Here, a planner manages multiple inventories to meet stochastic demand while minimizing costs that are coupled through nonlinear ordering functions reflecting economies of scale. Our main contributions are: (i) a surprising equivalence result showing that optimal stationary base-stock levels for individual locations remain unchanged despite the coupling when restricting attention to decoupled strategies; (ii) tight performance bounds for simple decoupled policies relative to optimal coupled policies, revealing that the worst-case ratio depends primarily on the degree of nonlinearity in the cost function and scales with the number of locations for systems with fixed costs; and (iii) analysis of practical online algorithms that achieve competitive performance without solving complex dynamic programs. Numerical simulations demonstrate that while decoupled policies significantly outperform their worst-case guarantees in typical scenarios, they still exhibit meaningful suboptimality compared to fully coordinated strategies. These results provide actionable guidance for system operators navigating the trade-off between control complexity and operational efficiency in coupled systems.
\end{abstract}

\section{Introduction}
Industrial refrigeration poses significant control challenges, particularly in cold storage facilities where multiple rooms must be maintained below critical temperature thresholds to preserve perishable goods. These systems use evaporators in each room that absorb heat from storage spaces and a centralized set of compressors—the dominant energy consumers—to meet the aggregate cooling load. Conventional control approaches rely on simple temperature setpoints that determine whether individual evaporators are on or off, creating a decoupled system where evaporator operation depends solely on conditions within each room \cite{hovgaard2012model, wisniewski2014analysis, mota2015commercial}. However, cooling multiple rooms simultaneously could be more energy efficient because compressors operate most efficiently at maximum capacity \cite{manske2001evaporative, reindl2013sequencing}. Although technologies now exist to enable such coordinated control strategies \cite{CrossnoKaye_2025}, it is crucial to determine whether the resulting efficiency gains would justify the implementation costs.

Cold storage facilities can be thought of as \emph{multi-location inventory control systems} \cite{veinott1965optimal, balintfy1964basic, federgruen1984approximations}, a model from operations research where a planner manages multiple inventories of an item that is depleted by a random demand process over time; the core challenge is to balance the competing costs of holding excess inventory, incurring backlogs when demand exceeds the current inventory level, and placing orders for additional inventory. In this framework, cooling capacity corresponds to the `inventory level' and heat loads represent `demands' that deplete this inventory. The `thermal battery' concept, i.e., intentionally overcooling products below required thresholds, parallels holding excess stock in traditional inventory systems, insufficient cooling is akin to backlogs, and the energy consumed by cooling operations represents ordering. When multi-location inventory systems are \emph{decoupled} with independent demands and costs, they reduce to separate single-location systems, and decoupled control policies (like the setpoint control approaches used in refrigeration systems) are often effective.
% In many cases, these single-location problem instances are known to have simple optimal policies (such as so-called base-stock or $(s, S)$ policies), analogous to the setpoint control approaches used in refrigeration systems.
However, when the systems are coupled through their costs, as in refrigeration systems where compressor efficiency creates interdependencies, the problem and its optimal solutions may become more complex.

Inventory control provides a well-established framework for studying coupling in these kinds of systems.
% In its classic formulation, a planner manages the inventory of a single item that is depleted by a random demand process over time. The core challenge is balancing the competing costs of holding excess inventory, incurring stockouts or backlogs when demand exceeds available supply, and placing orders.
Decades of research on the problem have focused on characterizing optimal policies under various assumptions on the demand process and cost functions for both single-location and multi-location systems. The key insight from this literature is that for typical single-location problems, simple control policies such as base-stock policies (order up to a specific threshold $S$ if inventory falls below $S$) or $(s, S)$ policies (order up to $S$ if inventory falls below $s$) are provably optimal \cite{arrow1951optimal, scarf1960optimality}.

The multi-location inventory problem generalizes the standard problem to multiple locations. When the systems are decoupled, one can achieve optimal or near-optimal performance with decoupled policies \cite{arrow1951optimal, veinott1965optimal}. However, in coupled systems where the cost in one location depends on another, such as when per-unit costs decrease with larger order quantities \cite{porteus1971optimality, benjaafar2018optimal}, the problem is no longer decomposable, and simple decoupled policies may be significantly suboptimal. Several works have characterized optimal \cite{balintfy1964basic, veinott1965optimal, eppen1981centralized, erkip1990optimal,liu2012decision} or near-optimal \cite{hochstaedter1970approximation, federgruen1984approximations,federgruen1984computational} policies for various kinds of coupled multi-location systems. Importantly, these policies can be complex and difficult to implement in practice.

The focus of this work is to study how well simple, decoupled policies perform in coupled multi-location inventory systems like the refrigeration example above. In particular, we focus on systems where demand and holding/backlog costs are independent across locations, but ordering costs are a nonlinear function of the cumulative order quantity across all locations in each period. Our results are summarized as follows:

(i) Our first contribution given in Theorem~\ref{thm:stationary_base_stock} demonstrates 
% a surprising result about decoupled base-stock policies in multi-location inventory control problems.
% In this context, a base-stock policy refers to a decoupled control strategy where each location independently follows its own base-stock rule, ordering up to a predetermined threshold whenever inventory falls below that level.
% We prove
that the optimal stationary base-stock policy for a particular location in a multi-location inventory control problem is equivalent to the optimal stationary base-stock policy for that location in isolation. In other words, despite the coupling through nonlinear ordering costs, the existence of other locations does not influence the optimal base-stock level for a given location. This result significantly simplifies the process of finding optimal stationary base-stock policies, as it allows us to solve separate single-location problems rather than a joint multi-location problem.

(ii) Our second contribution given in Theorem~\ref{thm:main} provides tight bounds on the performance gap between the optimal policy and the optimal decoupled base-stock and $(s, S)$ policies in multi-location inventory control problems. 
% We characterize the competitive ratio, i.e., the worst-case ratio between the cost of the optimal base-stock policy and the cost of the optimal policy, i.e., for several important classes of nonlinear ordering cost functions.
% We characterize the worst-case ratio between the cost of the optimal decoupled base-stock or $(s, S)$ policy and the cost of the optimal policy for several important classes of nonlinear ordering cost functions.
Specifically, for sector-bounded ordering costs, we show that this ratio depends only on the degree of nonlinearity of the ordering cost function, while for affine-bounded ordering costs, it also scales with the number of locations. These bounds precisely quantify the trade-off between policy simplicity and performance across different system structures, providing guidance on when decoupled base-stock policies are nearly optimal and when more complex coupled policies may be necessary.

(iii) Our third contribution extends our analysis to a simple online policy as an alternative to optimal decoupled base-stock or $(s, S)$ policies, which are structurally simple but can be computationally difficult to identify, especially in the case where demands are correlated over time~\cite{dong2003optimal}. In Proposition~\ref{prop:online_algs}, we analyze the randomized cost-balancing policy~\cite{levi2013approximation}, an efficiently computable heuristic that balances ordering, holding, and backlog costs. We prove that this algorithm achieves competitive ratios that scale predictably with our bounds for optimal decoupled policies, thereby demonstrating that effective inventory control is possible even with online decision rules. These results are particularly valuable for practical applications where computational resources are limited or real-time decisions are required.

We validate these theoretical results through numerical simulations which reveal that while simple decoupled policies and online algorithms outperform their worst-case bounds in practice, they still exhibit significant suboptimality compared to optimal coupled policies. Although our motivation stems from industrial refrigeration systems, we recognize that real refrigeration dynamics involve additional complexities not captured in our model. However, our purpose is not to accurately model refrigeration systems but rather to provide a first-order understanding of how decoupled control policies perform in inherently coupled environments. The fundamental insight of this paper is a precise characterization of the performance guarantees achievable with structurally simple control policies like setpoint control.
% , such as the setpoint control commonly used in refrigeration systems, 
These findings provide practical guidance for system operators facing the trade-off between control complexity and operational efficiency across various domains with similar coupling structures.

\begin{comment}
The remainder of this paper is organized as follows. Section \ref{sec:model} formalizes the model of the multi-location inventory problem. Section \ref{sec:policy_comparison} presents our main theoretical results on constant-factor bounds for decoupled policies, and Section \ref{sec:online_algs} extends these results to online algorithms. Finally, in Section \ref{sec:sims}, we evaluate the performance of the decoupled policies and online algorithms in numerical simulations.
\end{comment}

\section{Multi-Location Inventory Problem}\label{sec:model}

In the multi-location inventory problem, one seeks to meet stochastic demand while balancing ordering costs, holding costs, and backlog costs across multiple locations. The system is described by the $M$-dimensional linear dynamics
\begin{equation}\label{eq:dynamics}
    x_{k + 1} = x_k + u_k - w_k,
\end{equation}
where $x_k \in \mathbb{R}^M$ is the inventory level at each location at stage $k$, $u_k \in \mathbb{R}_{\geq 0}^M$ is the control input (ordering), and $w_k$ is an $M$-dimensional random variable that represents stochastic demand. We assume that each $w_k^i$ is finite-mean, independent, nonnegative, and has bounded support, where we use the superscript $i \in \{1, ..., M\}$ to index the relevant quantity for the $i$\textsuperscript{th} location.

At each stage, we face two competing costs: the cost of ordering, which we represent as $c : \mathbb{R}_{\geq 0} \to \mathbb{R}_{\geq 0}$, and the cost of holding inventory or experiencing backlogged demand, which we represent as $r :  \mathbb{R}^M \to \mathbb{R}_{\geq 0}$. We assume that $c$ is a function of the total order quantity across all the individual locations\footnote{We assume that the cost is an unweighted sum of the order amount at each location, but all of our results extend to the setting where the ordering cost is $c(\sum_{i=1}^M \alpha^i z^i)$ for some $\alpha > 0$.}, i.e., $c(z)$ is a shorthand notation for $c(\sum_{i=1}^M z^i)$, and that $c$ is locally bounded and lower semicontinuous.
% \james{It could be helpful to label the assumptions we use everywhere as Assumption 1 etc, then we can add `under Assumption 1...' to all the theorem statements to make them more precise.}
Additionally, we assume that the holding/backlog costs are additively separable, i.e., $r(z) = \sum_{i = 1}^M r^i(z^i)$, where $r^i: \mathbb{R} \to \mathbb{R}_{\geq 0}$ is locally bounded, radially unbounded, and convex. The average costs over an $N$-period horizon are given by
\begin{equation}\label{eq:multi_inv_prob}
    \mathbb{E} \left[ \frac{1}{N} \sum_{k=0}^{N-1} \left( c(u_k) + r(x_k + u_k - w_k) \right) \right].
\end{equation}
For the infinite-horizon setting, we follow the approach of~\cite{schal1993average} and average costs by taking the limit superior of \eqref{eq:multi_inv_prob} as $N \to \infty$. We adopt the notation $P = (c, r, W)$ to completely describe a multi-location inventory problem, where the dimension of $W$, an $M$-by-$N$ random variable whose $(i, k)$-th element is $w_{k - 1}^i$, indicates the number of inventories and the length of the horizon, and we use $\mathcal{P}^M$ to denote the set of all $M$-location inventory problems with finite or infinite horizon that meet the above assumptions on $c$, $r$, and $W$.

We define a \emph{policy} $\pi = \{ \mu_k \}_{k \, \geq\, 0}$ as a sequence of functions $\mu_k : \mathbb{R}^M \to \Delta(\mathbb{R}^M_{\geq 0})$ that map a state $x_k$ to a probability distribution over feasible controls $u_k^i \geq 0$ for all $i$; if $\mu_k = \mu$ for all $k$, then $\pi$ is said to be \emph{stationary}. The expected cost for a given policy $\pi$ from initial condition $x_0$ is given by $J_\pi(x_0 \given P)$, which is defined as
\begin{equation}\label{eq:multi_inv_prob_w_pol}
    \E \left[ \frac{1}{N} \sum_{k = 0}^{N-1} \left( c(\mu_k(x_k)) + r(x_k + \mu_k(x_k) - w_k) \right) \right]
\end{equation}
with appropriate modifications for the infinite-horizon setting. An \emph{optimal} policy $\pi_*$ for a given problem $P$ is one which incurs expected cost less than or equal to all other policies for all initial conditions, i.e.,
\begin{equation}\label{eq:optimality_condition}
    \pi_* \in \argmin_{\pi \in \Pi^M} J_{\pi}(x_0 \given P), \quad \forall \ x_0 \in \mathbb{R}^M
\end{equation}
where $\Pi^M$ is the set of all $M$-dimensional policies; our assumptions on $c$, $r$, and $W$ guarantee the existence of $\pi_*$ for both finite~\cite[Ch.~3]{bertsekas1996stochastic} and infinite-horizon problems \cite{schal1993average}.
% 
% \james{Invoke standard MPD results for 1) the existence of $\pi_*$, and 2) restricting our attention to deterministic policies?}
% 
A policy $\pi$ is said to be \emph{decoupled} if $\mu_k$ is component-wise separable for all $k$, i.e., $\mu_k(z_k) = \{\mu_k^1(z_k^1), ..., \mu_k^M(z_k^M) \}$; otherwise, it is \emph{coupled}. We occasionally abuse notation and write a decoupled policy as $\pi = \{\pi^1, ..., \pi^M \}$, where $\pi^i = \{ \mu_k^i \}_{k \geq 0}$.

\subsection{Motivating Examples}
\begin{figure*}
    \centering
    \includegraphics[width=\textwidth]{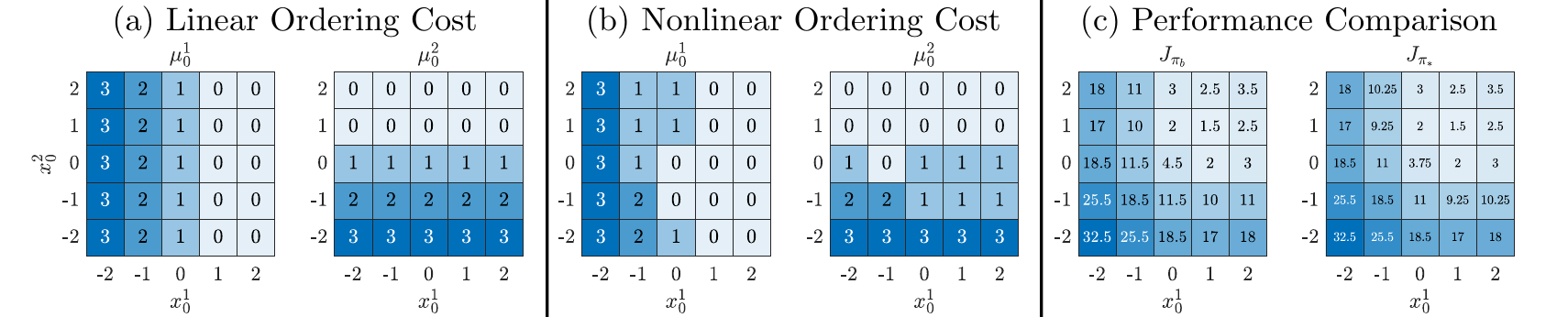}
    \caption{We consider $M = 2$ inventories, integer-valued parameters, a finite horizon $N = 2$, holding/backlog cost $r(z) = \max \{0,z\} + 10 \max \{0,-z\}$, and stochastic demand $w_k$ that takes values zero or one with equal probability. (a) shows the optimal policy for the linear ordering cost $c(z) = 2z$. (b) shows the optimal policy for the nonlinear ordering cost function~\eqref{eq:c_nonlinear}. (c) shows the expected cost from each initial condition for the nonlinear ordering cost function~\eqref{eq:c_nonlinear} under the policies in (a),(b).}
    \label{fig:num_ex}
\end{figure*}

Fig.~\ref{fig:num_ex} contains some numerical examples that illustrate optimal policies under linear and nonlinear ordering costs and motivate our study of the suboptimality of decoupled policies. We define \emph{linear} ordering costs as functions of the form $c(z) = mz$ where $m \geq 0$ is the per-unit cost. 

Fig.~\ref{fig:num_ex}(a) shows the optimal policy for a 2-location problem with ordering cost $c(z) = 2z$ and has two notable features.  First, because the multi-location inventory problem~\eqref{eq:multi_inv_prob} is coupled only through the ordering cost function, for linear ordering costs the problem decouples into $M$ independent single-location problems. As a result, the optimal policy in this case is also decoupled. Second, it is well-known in the inventory control literature that for single-location problems with linear ordering costs, so-called base-stock policies are optimal \cite{arrow1951optimal}.
\begin{defn}
    A policy $\pi \in \Pi^1$ is a \emph{base-stock policy} if
    \begin{equation}\label{eq:base_stock_defn}
        \mu_k(x_k) = \max\{ S_k - x_k, 0 \}
    \end{equation}
    for some $S_k \in \mathbb{R}$ for all $k \geq 0$.
\end{defn}
% Base-stock policies are characterized by a single scalar, the order-up-to level $S_k$, at each stage $k$. The policy consists entirely of ordering up to $S_k$ in states $x_k \leq S_k$ and ordering nothing in the remaining states.
\noindent Fig.~\ref{fig:num_ex}(a) shows that the optimal policy for each location at time $k = 0$ is a base-stock policy where $S_0 = 1$. 

Fig.~\ref{fig:num_ex}(b) shows that the optimal policy for a nonlinear ordering cost, in this case
\begin{equation}\label{eq:c_nonlinear}
    c(z) = \begin{cases}
        2z, \quad & z \in \{0,1\} \\
        4z, \quad & z \in \{2,3,4\},
    \end{cases}
\end{equation}
can be coupled and have quite general structure. Note that in the $x_0^1 = x_0^2 = 0$ state it is optimal to order 1 unit in either location (but not both) so the two locations are indistinguishable, i.e., there is symmetry about the anti-diagonal.  For this particular nonlinear ordering cost function~\eqref{eq:c_nonlinear}, the best (decoupled) base-stock policy also happens to be the one shown in Fig.~\ref{fig:num_ex}(a). 

Fig.~\ref{fig:num_ex}(c) shows the increase in expected cost associated with using the policy from Fig.~\ref{fig:num_ex}(a) instead of the optimal policy from Fig.~\ref{fig:num_ex}(b). To characterize the suboptimality in cases like this, our main results provide bounds on the worst-case cost increase of using decoupled single-location policies for multi-location problems for broad classes of nonlinear ordering cost functions.

\section{Policy Comparison}\label{sec:policy_comparison}
In this section, we compare the performance of optimal and decoupled policies in the multi-location inventory problem \eqref{eq:multi_inv_prob_w_pol} for various ordering cost functions. Our results in this section focus specifically on the infinite horizon setting, i.e., $N = \infty$.

\subsection{Stationary Base-Stock Policies}
The motivating examples in the previous section highlight that complex, coupled policies may perform better than simple, decoupled ones. In this section, we momentarily disregard questions of policy complexity and ask instead about the effects of coupling within an extremely simple class of policies.
% , as the state space of the corresponding dynamic program grows exponentially in $M$.
% \james{This is only true if we are discretizing the true state space $\real^M$ -- this seems to be implied, but have we ever explicitly stated it?}
% Typically, this simplicity comes at the expense of optimality, but in this section, we establish that for certain classes of policies, this is not the case.
Consider the class of \emph{decoupled stationary base-stock policies} $\Pi_\textup{Y}^{M} \subset \Pi^M$, where a decoupled policy $\pi = \{\pi^1, ..., \pi^M \}$ is said to be a decoupled stationary base-stock policy if $\pi^i$ is a stationary base-stock policy for all $i$. Our first result demonstrates that when optimizing within this class for infinite horizon problems, coupling can be neglected, i.e., there is no difference between solving for optimal base-stock policies jointly over $\Pi_\textup{Y}^{M}$ or individually over $\Pi_\textup{Y}^{1}$.

\begin{theorem}[Equivalence of stationary base-stock policies]\label{thm:stationary_base_stock}
    Consider an infinite horizon $M$-location inventory problem $P = (c, r, W) \in \mathcal{P}^M$ with initial condition $x_0 \in \mathbb{R}^M$. The jointly-optimized stationary base-stock policy
    \begin{equation}\label{eq:jointly_optimized}
        \pi_\star \in \argmin_{\pi \in \Pi_\textup{Y}^M} J_{\pi}(x_0 \given P)
    \end{equation}
    and the individually-optimized stationary base-stock policy $\pi_\circ = \{\pi_\circ^1, ..., \pi_\circ^M\}$, where $P^i = (c, r^i, W^i)$ and
    \begin{equation}\label{eq:individually_optimized}
        \pi_\circ^i \in \argmin_{\pi \in \Pi_\textup{Y}^1} J_{\pi}(x_0^i \given P^i),
    \end{equation}
    exist and satisfy $J_{\pi_\circ}(x_0 \given P) = J_{\pi_\star}(x_0 \given P)$.
\end{theorem}

% 
% \james{It's a bit surprising that we don't need to change the cost function when constructing the individual $P^i$ from the global $P$; some discussion would be appreciated.}
% \james{State conditions more clearly -- are we still assuming sector-bounded?}
% 
% \james{I think we intend to mean that $\Pi_B$ are whole-system policies formed from a product of independent base-stock policies for each location. However, it is not totally clear here, because we have not explicitly stated that our concept of base-stock policies doesn't include some notion of a base-stock policy that applies to the entire coupled system.}
% 

Theorem \ref{thm:stationary_base_stock} highlights a key property of stationary base-stock policies in infinite-horizon problems: within this class, solving for policies independently is optimal. Importantly, this means that the coupling in the ordering cost can and should be disregarded, as one cannot achieve better performance by taking it into consideration. Furthermore, this result is quite general, as it holds for any locally bounded, nonnegative, and lower semicontinuous ordering cost function.

\begin{proof}
    Let $S^i$ denote the base-stock level in inventory $i$ defined by $\pi_\star$, and assume without loss of generality that $S^i \geq 0$ for all $i$. Given that the support of the demand $w_k^i$ is nonnegative, we can write the objective function in \eqref{eq:jointly_optimized} with the initial condition $x_0^i = 0$ for all $i$ as
    % \begin{align*}
    %     \min_{S \in \mathbb{R}^M} \limsup_{N \to \infty} \frac{1}{N} \mathbb{E} \Bigg[ \sum_{k=0}^N \Bigg( c\left( \sum_{i = 1}^M (S^i - x_k^i) \right) \\
    %     + \sum_{i = 1}^M r^i(S^i - w_k^i) \Bigg) \Bigg].
    % \end{align*}
    \begin{equation*}
        \limsup_{N \to \infty} \frac{1}{N} \mathbb{E} \Bigg[ \sum_{k=0}^N \Bigg( c\left( \sum_{i = 1}^M (S^i - x_k^i) \right)
        + \sum_{i = 1}^M r^i(S^i - w_k^i) \Bigg) \Bigg].
    \end{equation*}
    From \eqref{eq:dynamics} and \eqref{eq:base_stock_defn}, we have that for all $k \geq 1$,
    \begin{align*}
        u_k^i = S^i - x_k^i &= S^i - (x_{k - 1}^i + u_{k - 1}^i - w_{k - 1}^i ) \\
        &= S^i - (x_{k - 1}^i + S^i - x_{k - 1}^i - w_{k - 1}^i ) = w_{k - 1}^i,
    \end{align*}
    so we can separate the $k = 0$ terms and rewrite the optimization problem as
    \begin{align*}
        \min_{S \in \mathbb{R}^M} \limsup_{N \to \infty} \frac{1}{N} \mathbb{E} \Bigg[  c\left( \sum_{i = 1}^M (S^i - x_0^i) \right) + \sum_{i = 1}^M r^i(S^i - w_0^i) \\
        + \sum_{k=1}^N \left( c\left( \sum_{i = 1}^M w_{k - 1}^i \right) + \sum_{i = 1}^M r^i(S^i - w_k^i) \right) \Bigg].
    \end{align*}
    The expected ordering cost terms are a fixed constant for all $k \geq 1$, and since $c$ and $r$ are locally bounded, the total cost incurred in the first period is finite. Hence, we can equivalently rewrite the problem as
    \begin{align*}
        \min_{S \in \mathbb{R}^M} \limsup_{N \to \infty} \frac{1}{N} \mathbb{E} \left[ \sum_{k=1}^N  \sum_{i = 1}^M r^i(S^i - w_k^i) \right].
    \end{align*}
    This problem can be decomposed into $M$ independent single-location inventory control problems, and the same analysis can be applied to verify that the single-inventory problem corresponding to location $i$ is precisely problem \eqref{eq:individually_optimized}, meaning that the expected costs of $\pi_\star$ and $\pi_\circ$ from the zero initial condition are identical. The existence of both $\pi_\star$ and $\pi_\circ^i$ follows from the fact that the limit superior of the objective exists and is a nonnegative, radially unbounded, and proper function in $S$, so it has an infimum; furthermore, the infimum can be achieved because it is also lower semicontinuous. 
    
    We extend the analysis to nonzero initial conditions as follows. For any initial condition where $x_0^i < 0$, a stationary threshold policy orders up to $S^i$ in one period, so by the same argument as above, these costs are finite and negligible. For any initial condition where $x_0^i > 0$, we consider two cases: if $\sum_{k = 1}^\infty \mathbb{E}[w_k^i] = \infty$, then $x_k^i \leq S^i$ for some $k < \infty$ with probability 1, and the finite costs associated with these periods can be eliminated in a similar fashion. Otherwise, if $\sum_{k = 1}^\infty \mathbb{E}[w_k^i] < \infty$, $S^i = 0$ is trivially optimal for both the individually- and jointly-optimized problems.
\end{proof}

\subsection{Base-Stock and $(s, S)$ Policies}\label{sec:main}

Unfortunately, the results of the previous section are specific to stationary base-stock policies in infinite horizon problems. It is not difficult to show that when one allows for more general policies or when one considers a finite horizon, coupling cannot always be neglected without some performance degradation. Hence, in this section, we study both complexity and coupling as we quantify the performance gap between optimal coupled and simple decoupled policies. We focus specifically on the classes of \emph{decoupled base-stock} and \emph{decoupled $(s, S)$} policies $\Pi_\textup{B}^{M} \subset \Pi^{M}$ and $\Pi_\textup{S}^{M} \subset \Pi^{M}$, respectively; a policy $\pi = \{ \pi^1, ..., \pi^M \}$ is a decoupled base-stock (resp., $(s, S)$) policy if $\pi^i$ is a base-stock (resp., $(s, S)$) policy for all $i$.
\begin{defn}
    A policy $\pi \in \Pi^1$ is an \emph{$(s, S)$ policy} if
    \begin{equation}
        \mu_k(x_k) = \begin{cases}
            S_k - x_k & x_k < s_k, \\
            0 & x_k \geq s_k
        \end{cases}
    \end{equation}
    for some $s_k, S_k \in \mathbb{R}$ where $s_k \leq S_k$ for all $k \geq 0$.
\end{defn}
\noindent Importantly, $(s, S)$ policies are known to be optimal for single-location inventory problems with both fixed and per-unit ordering costs, i.e., $c(z) = K \1(z) + mz$ where $\1(z) = 1$ if $z > 0$ and $0$ otherwise \cite{scarf1960optimality}. We refer to such ordering cost functions as \emph{affine}.

The results in this section focus on two broad classes of nonlinear ordering cost functions:

\begin{defn}
    A function $c : \mathbb{R}_{\geq 0} \to \mathbb{R}_{\geq 0}$ is \emph{sector-bounded} if $lz \leq c(z) \leq hz$
    for some $h \geq l > 0$. 
\end{defn}

\begin{defn}
    A function $c : \mathbb{R}_{\geq 0} \to \mathbb{R}_{\geq 0}$ is \emph{affine-bounded} if $
    K_l \1(z) + lz \leq c(z) \leq K_h \1(z) + hz$
    for some $h \geq l > 0$ and $K_h \geq K_l > 0$.
\end{defn}

Our main result establishes bounds on the worst-case performance of decoupled base-stock policies for infinite horizon problems\footnote{Note that in the finite horizon setting, \eqref{eq:affine_bounded_bound} still holds; however, \eqref{eq:sector_bounded_bound} holds as an upper bound rather than as an equality.} with sector-bounded ordering costs and decoupled $(s, S)$ policies for problems with affine-bounded ordering costs.

\begin{theorem}[Suboptimality of base-stock and $(s,S)$ policies]\label{thm:main}
    \leavevmode \vspace{-\baselineskip} 
    \begin{enumerate}[wide = 0pt]
        \item Let $\mathcal{P}_\textup{SB}^M \subset \mathcal{P}^M$ denote the set of all $M$-location inventory problems where the ordering cost is sector-bounded with $h \geq l > 0$. We have that
        \begin{equation}\label{eq:sector_bounded_bound}
            \sup_{x_0 \in \mathbb{R}^M,P \in \mathcal{P}_\textup{SB}^M} \min_{\pi \in \Pi_\textup{B}^M} \frac{J_{\pi}(x_0 \given P)}{J_{\pi_*}(x_0 \given P)} = \frac{h}{l}.
        \end{equation}
        \item Let $\mathcal{P}_\textup{AB}^M \subset \mathcal{P}^M$ denote the set of all $M$-location inventory problems where the ordering cost is affine-bounded with $h \geq l > 0$ and $K_h \geq K_l > 0$. We have that \begin{equation}\label{eq:affine_bounded_bound}
         \sup_{x_0 \in \mathbb{R}^M,P \in \mathcal{P}_\textup{AB}^M} \min_{\pi \in \Pi_\textup{S}^{M}} \frac{J_{\pi}(x_0 \given P)}{J_{\pi_*}(x_0 \given P)} \leq M {\max} \left\{ \frac{K_h}{K_l} {,} 
         \frac{h}{l} \right\}.
    \end{equation}
    \end{enumerate}
\end{theorem}

Theorem \ref{thm:main} establishes several interesting properties of optimal decoupled base-stock policies for problems with sector-bounded ordering costs. First, observe that the ratio \eqref{eq:sector_bounded_bound} grows with the ratio $h/l$, indicating that optimal decoupled policies may perform worse with `more nonlinear' functions; note that when the ordering cost is linear, $h/l = 1$. Second, the proof is constructive, meaning that one can identify a decoupled base-stock policy that achieves this bound by solving only $M$ single-location problems with identical linear ordering costs. Third, the tightness of the result means that no base-stock policy achieves a ratio strictly less than $h/l$ for all problem instances, implying that this policy is not only reasonable to implement but also optimal (with respect to worst-case performance) among base-stock policies.
% 
% \james{ Are you sure about this explanation? I speculate one could construct a $P$ with a tight $h,l$ sector bound such that $J_{\pi_*} = J_\pi$, for example by making $c$ linear except on a set of measure zero, and choosing $W$ appropriately. To me, the tightness only asserts exactly the claim from the theorem: \emph{there exist} problems where no base-stock policy is better than $h/l$, so this is the best \emph{single} upper bound that holds across \emph{all} problems. That does not mean problem instances where decoupled can do better than $h/l$ \emph{do not exist}.}
% 

Theorem \ref{thm:main} also establishes important properties of optimal decoupled $(s, S)$ policies for problems with affine-bounded ordering costs. Here, the bound \eqref{eq:affine_bounded_bound} also grows with the ratios $h/l$ and $K_h / K_l$, indicating that the degree of nonlinearity affects worst-case performance similarly, and the proof is also constructive. However, the key distinction between the two is that the upper bound in \eqref{eq:affine_bounded_bound} scales with the number of locations. This suggests that for systems with several locations and fixed costs—such as cold storage facilities, where there may be several rooms and compressors typically have minimum operating thresholds \cite{reindl2013sequencing}—decoupled policies may perform especially poorly\footnote{That said, we acknowledge that this bound is not tight, and we leave it to future work to show that this ratio is also lower-bounded by an increasing function of $M$.}.
% 
% \james{Would be nice to recall the cooling systems example here.}
% 

\begin{proof}
    We first establish \eqref{eq:sector_bounded_bound}.
    We proceed constructively. Taking $c_l(z) = lz$, let $\pi_\Hsquare$ denote the (optimal) decoupled base-stock policy for the multi-location problem $P_l = (c_l, r, W)$. We have that
    \begin{align*}
        \frac{J_{\pi_\Hsquare}(x_0 \given P)}{J_{\pi_*}(x_0 \given P)}
        % \leq \frac{J_{\pi_\Hsquare}(x_0 \given P)}{J_{\pi_*}(x_0 \given P_l)}
        \hintedrel[l_star]{\leq} \frac{J_{\pi_\Hsquare}(x_0 \given P)}{J_{\pi_\Hsquare}(x_0 \given P_l)}
        \leq \frac{J_{\pi_\Hsquare}(x_0 \given P_h)}{J_{\pi_\Hsquare}(x_0 \given P_l)}
        \hintedrel[sec_bd]{\leq} \frac{h}{l}
    \end{align*}
    where $P_h = (c_h, r, W)$ with $c_h(z) = hz$. Here, (\hintref{l_star}) holds because the costs under $P$ are greater than or equal to the costs under $P_l$, and (\hintref{sec_bd}) follows from the fact that since the same policy $\pi_\Hsquare$ is enacted on both $P_l$ and $P_h$, the only differences in their overall costs must result from differences in their ordering costs which are bounded by the sector. Because these inequalities hold for any $x_0$ and $P$ with sector-bounded ordering cost, they hold for the supremum over $x_0$ and $P$. To establish equality, we provide an example of a generic infinite horizon $M$-location inventory problem in Appendix~\ref{app:tightness_thm_1} such that for any $\varepsilon > 0$, the expected cost of any base-stock policy is at least $\frac{h}{l + \varepsilon}$ greater than the expected cost of the optimal policy for some initial condition. Taking the limit as $\varepsilon \to 0$ yields that the supremum exists and is $\frac{h}{l}$.

    To show \eqref{eq:affine_bounded_bound}, we first consider the setting where $K_l = K_h$ and $l = h$ (i.e., the function is affine). Let $\pi_\diamond$ denote the decoupled $(s, S)$ policy $\{\pi_\diamond^1, ..., \pi_\diamond^M\}$ where $\pi_\diamond^i$ is optimal for $P^i = (c, r^i, W^i)$. Let $j \in \argmax_{i \in \{1, ..., M\}} \ J_{\pi_\diamond^i}(x_0 \given P^i)$ denote the single location that incurs the maximum expected cost under $\pi_\diamond$ for the initial condition $x_0$. Then, we have that
    \begin{equation}\label{eq:affine_bound}
        \frac{J_{\pi_\diamond}(x_0 \given P)}{J_{\pi_*}(x_0 \given P)}
        \hintedrel[single_multi]{\leq} \frac{J_{\pi_\diamond}(x_0 \given P)}{J_{\pi_\diamond^j}(x_0 \given P^j)}
        \leq \frac{M J_{\pi_\diamond^j}(x_0 \given P^j)}{J_{\pi_\diamond^j}(x_0 \given P^j)}
        = M
    \end{equation}
    where we have (\hintref{single_multi}) because adding additional locations cannot decrease the expected cost. 

    We now consider the general setting. Define $c_l(z) = K_l \1(z) + lz$ and $c_h(z) = K_h \1(z) + hz$, and let $\pi_\diamond$ instead denote the decoupled policy $\{\pi_\diamond^1, ..., \pi_\diamond^M\}$ where $\pi_\diamond^i$ is optimal for $P^i = (c_h, r^i, W^i)$. Let $\pi_l$, $\pi_*$, and  $\pi_h$ denote the optimal policies for the problems $P_l = (c_l, r, W)$, $P$, and $P_h = (c_h, r, W)$, respectively. First, we apply some trivial bounds on the ratio of the expected costs between $\pi_\diamond$ and $\pi_*$:
    \begin{align*}
        \frac{J_{\pi_\diamond}(x_0 \given P)}{J_{\pi_*}(x_0 \given P)}
        % \leq \frac{\frac{1}{M} J_{\pi_\diamond}(x_0 \given P)}{J_{\pi_*}(x_0 \given  P_l)}
        \leq \frac{J_{\pi_\diamond}(x_0 \given P)}{J_{\pi_l}(x_0 \given P_l)} \leq \frac{J_{\pi_\diamond}(x_0 \given P_h)}{J_{\pi_l}(x_0 \given P_l)}.
    \end{align*}
    We further bound the last expression by
    \begin{align*}
        \frac{J_{\pi_\diamond}(x_0 \given P_h)}{J_{\pi_l}(x_0 \given P_l)}
        \hintedrel[aff_bd]{\leq}
        \frac{M J_{\pi_h}(x_0 \given P_h)}{J_{\pi_l}(x_0 \given P_l)}
        \hintedrel[pol_chg]{\leq}
        \frac{M J_{\pi_l}(x_0 \given P_h)}{J_{\pi_l}(x_0 \given P_l)}
    \end{align*}
    where (\hintref{aff_bd}) follows from \eqref{eq:affine_bound} and (\hintref{pol_chg}) follows from the fact that $\pi_l$ is suboptimal for $P_h$. All that remains is to bound the difference in costs between the problems $P_l$ and $P_h$ when $\pi_l$ is applied. Consider performing the cost transformation described in Appendix~\ref{app:cost_transformation}
    % 
    % \james{ref. specific section}
    % 
    on $P_l$ and $P_h$. We can use~\eqref{eq:cost_transformation} to write
    \begin{equation*}
        \frac{J_{\pi_l}(x_0 \given P_h)}{J_{\pi_l}(x_0 \given P_l)} = \frac{J_{\pi_l}(x_0 \given \hat{P}_h) + h \gamma}{J_{\pi_l}(x_0 \given \hat{P}_l) + l \gamma} = \frac{\alpha + K_h \beta + h \gamma}{\alpha + K_l \beta + l \gamma}
    \end{equation*}
    where $\alpha \geq 0$ represents expected holding/backlog costs, $\beta = \frac{1}{N} \sum_{k=0}^{N-1}\1( \sum_{i=1}^M u_k^i) \geq 0$ represents fixed ordering costs, and $\gamma = \E\left[\frac{1}{N} \sum_{k=0}^{N-1} \sum_{i=1}^M w_k^i \right] \geq 0$ represents unit ordering costs (with appropriate limits for the infinite horizon setting). 
    We proceed by considering two cases: 
    \begin{enumerate}
        \item If $K_l h \leq K_h l$, then one can show that $\frac{\alpha + K_h \beta + h \gamma}{\alpha + K_l \beta + l \gamma} \leq \frac{K_h}{K_l}$.
        \item If $K_l h > K_h l$, then it can similarly be shown that $\frac{\alpha + K_h \beta + h \gamma}{\alpha + K_l \beta + l \gamma} < \frac{h}{l}$.
    \end{enumerate}
    % If $K_l h \leq K_h l$, then one can show that $\frac{\alpha + K_h \beta + h \gamma}{\alpha + K_l \beta + l \gamma} \leq \frac{K_h}{K_l}$. Otherwise, if $K_l h > K_h l$, then it can similarly be shown that $\frac{\alpha + K_h \beta + h \gamma}{\alpha + K_l \beta + l \gamma} < \frac{h}{l}$.
    % 1) $K_l h \leq K_h l$: It can be shown that this implies
    % \begin{equation*}
    %     \frac{\alpha + K_h \beta + h \gamma}{\alpha + K_l \beta + l \gamma} \leq \frac{K_h}{K_l}.
    % \end{equation*}
    % 2) $K_l h > K_h l$: It can be shown that this implies
    % \begin{equation*}
    %     \frac{\alpha + K_h \beta + h \gamma}{\alpha + K_l \beta + l \gamma} < \frac{h}{l}.
    % \end{equation*}
    Thus, in either case the ratio is less than $\max \left\{ K_h/K_l, h/l \right\}$\footnote{Note that for the case of $h = l = 0$, the upper bound is $MK_h/K_l$ by the same argument.}. Once again, all of the inequalities hold for any $x_0$ and $P$ with the appropriate ordering cost function. Therefore, they hold for the supremum over $x_0$ and $P$. 
\end{proof}

\section{Online Algorithm}\label{sec:online_algs}
Computing an optimal policy for an inventory control problem can be challenging even in the single-location setting. For example, in the case of demand structures with correlation over time, the problem of solving for optimal policies can suffer from the curse of dimensionality~\cite{dong2003optimal}. However, for single-location finite-horizon problems with linear or affine ordering costs, there exists an efficient online algorithm that achieves constant-factor competitive ratios relative to optimal policies, even with correlated demands~\cite{levi2007approximation, levi2013approximation}. Here the term \emph{online algorithm} refers to policies that do not require pre-computation and can be efficiently implemented online. In this section, we analyze the application of this algorithm to multi-location problems with nonlinear ordering costs.

\subsection{Randomized Cost-Balancing Policy}
In~\cite{levi2013approximation}, the authors propose the \emph{randomized cost-balancing policy} $\pi_\vartriangle$ which is an online algorithm for single-location inventory control problems with linear or affine ordering costs and holding/backlog costs of the form $r(x) = a(\max\{0,x\}) + b(\max\{0,-x\})$ where $a,b \geq 0$. We define $\mathcal{Q}^M \subset \mathcal{P}^M$ as the set of all $M$-location finite-horizon inventory problems with ordering costs assumptions as in Section~\ref{sec:model} and holding/backlog costs as above. Under $\pi_\vartriangle$, at each time step one employs the following probabilistic ordering rule to balance in expectation the holding, backlog, and fixed cost $K$:
\begin{equation}
    \mu_\vartriangle = \begin{cases}
        \hat{u}_k, \quad & \theta_k \geq K,\ \textup{w.p.} \ 1 \\
        \tilde{u}_k, \quad & \theta_k < K,\ \textup{w.p.}\ p_k \\
        0, \quad & \theta_k < K,\ \textup{w.p.} \ 1-p_k
    \end{cases}
\end{equation}
where we suppress the dependency on $x_k$ and the remaining parameters depend on $H_k$ and $B_k$, proxy functions for holding and backlog costs that will be defined subsequently. The balancing quantity $\hat{u}_k$ yields balancing cost $\theta _k = \E [H_k(\hat{u}_k)] = \E [B_k(\hat{u}_k)]$. The holding-cost-K quantity $\tilde{u}_k$ sets $\E [H_k(\tilde{u}_k)] = K$. The probability $p_k$ solves
\begin{equation*}
    p_k K = p_k \E [B_k(\tilde{u}_k)] + (1-p_k)\E [B_k(0)].
\end{equation*}
Note that in the case of linear ordering costs, i.e., $K = 0$, the policy is deterministic. The rationale for the choice of $H_k$ and $B_k$ is that the holding costs associated with the order in a particular time step must be evaluated over the entire remaining horizon (over-ordering cannot be corrected) while the backlog costs correspond only to the next time step (under-ordering can be corrected). See~\cite{levi2013approximation} for details. At time step $k$, the holding cost over the remaining horizon for the units ordered now is 
\begin{equation*}
    H_k(u_k) = \sum_{n=k}^N a \max \{0,u_k - \max\{0, w_n - x_k\}\},
\end{equation*}
and the one-step backlog cost is 
\begin{equation*}
    B_k(u_k) = b (w_k - \max \{0, x_k + u_k\}).
\end{equation*}
% The rationale is that under-ordering can be remedied in the subsequent time step, but over-ordering cannot be corrected. 

% This policy balances the three types of costs. See~\cite{levi2013approximation} for details. 
The authors of~\cite{levi2013approximation} show that $\pi_\vartriangle$ is a tight 2-approximation for single-location inventory problems $\mathcal{Q}^1_L$ with linear ordering costs, i.e.,
\begin{equation}\label{eq:levi_1}
    \sup_{x_0 \in \mathbb{R},Q \in \mathcal{Q}^1_L} \frac{J_{\pi_\vartriangle}(x_0 \given Q)}{J_{\pi_*}(x_0 \given Q)} = 2,
\end{equation}
and a 3-approximation for single-location inventory problems $\mathcal{Q}^1_A$ with affine ordering costs, i.e.,
\begin{equation}\label{eq:levi_2}
    \sup_{x_0 \in \mathbb{R},Q \in \mathcal{Q}^1_A} \frac{J_{\pi_\vartriangle}(x_0 \given Q)}{J_{\pi_*}(x_0 \given Q)} \leq 3.
\end{equation}
For the $M$-location problem, we consider applying the randomized cost-balancing policy to each location independently, which we write as $\pi_\vartriangle = \{ \pi_\vartriangle^1, \dots, \pi_\vartriangle^M \}$. Combining~\eqref{eq:levi_1} and \eqref{eq:levi_2} with our results from Section~\ref{sec:main} yields the expected bounds on the suboptimality of $\pi_\vartriangle$:
\begin{prop}[Suboptimality of online algorithm]\label{prop:online_algs}
    \leavevmode % \vspace{-\baselineskip} 
    \begin{enumerate}[wide = 0pt]
        \item Let $\mathcal{Q}_\textup{SB}^M \subset \mathcal{Q}^M$ denote the set of all $M$-location inventory problems where the ordering cost is sector-bounded with $h \geq l > 0$. We have that
        \begin{equation}\label{eq:thm_online_1}
            \sup_{x_0 \in \mathbb{R}^M,Q \in \mathcal{Q}_\textup{SB}^M} \frac{J_{\pi_\vartriangle}(x_0 \given Q)}{J_{\pi_*}(x_0 \given Q)} \leq \frac{2h}{l}.
        \end{equation}
        \item Let $\mathcal{Q}_\textup{AB}^M \subset \mathcal{Q}^M$ denote the set of all $M$-location inventory problems where the ordering cost is affine-bounded with $h \geq l > 0$ and $K_h \geq K_l > 0$. We have that \begin{equation}\label{eq:thm_online_2}
         \sup_{x_0 \in \mathbb{R}^M, Q \in \mathcal{Q}_\textup{AB}^M} \frac{J_{\pi_\vartriangle}(x_0 \given Q)}{J_{\pi_*}(x_0 \given Q)} \leq 3M \max \left\{ \frac{K_h}{K_l}, 
         \frac{h}{l} \right\}.
    \end{equation}
    \end{enumerate}
\end{prop}
\begin{proof}
    We begin with~\eqref{eq:thm_online_1}. The proof is handled in two steps. First, we have the following bound when comparing $\pi_\vartriangle$ with $\pi_\Hsquare$, the optimal (decoupled) policy for the problem $Q_l = (c_l,r,W)$ with $c_l(z)=lz$:
    \begin{equation*}
        \frac{\frac{1}{2} J_{\pi_\vartriangle}(x_0 \given Q_h)}{J_{\pi_\Hsquare}(x_0 \given  Q_l)}
        \hintedrel[2_apprx]{\leq} \frac{J_{\pi_\Hsquare}(x_0 \given Q_h)}{J_{\pi_\Hsquare}(x_0 \given  Q_l)}
        \leq \frac{h}{l},
    \end{equation*}
    where $Q_h = (c_h,r,W)$ with $c_h(z) = hz$, and we use~\eqref{eq:levi_1} for (\hintref{2_apprx}). The second step is straightforward:
    \begin{equation*}
        \frac{\frac{1}{2} J_{\pi_\vartriangle}(x_0 \given Q)}{J_{\pi_*}(x_0 \given Q)}
        \leq \frac{\frac{1}{2} J_{\pi_\vartriangle}(x_0 \given Q)}{J_{\pi_\Hsquare}(x_0 \given  Q_l)}
        \leq \frac{\frac{1}{2} J_{\pi_\vartriangle}(x_0 \given Q_h)}{J_{\pi_\Hsquare}(x_0 \given  Q_l)}.
    \end{equation*}
    Now we move to showing~\eqref{eq:thm_online_2}. We will need the following bound for affine ordering costs that uses the decoupled $(s,S)$ policy $\pi_\diamond = \{\pi_\diamond^1,\dots,\pi_\diamond^M\}$ where $\pi_\diamond^i$ is optimal for $Q^i = (c, r^i, W^i)$ with $c(z) = K\1(z)+mz$:
        \begin{equation}\label{eq:online_affine}
        \frac{\frac{1}{3} J_{\pi_\vartriangle}(x_0 \given Q)}{J_{\pi_*}(x_0 \given Q)}
        \hintedrel[3_apprx]{\leq} \frac{J_{\pi_\diamond}(x_0 \given Q)}{J_{\pi_*}(x_0 \given Q)}
        \hintedrel[aff]{\leq} M.
    \end{equation}
    We use~\eqref{eq:levi_2} for (\hintref{3_apprx}) and (\hintref{aff}) follows from the proof of Theorem~\ref{thm:main}. Now we can proceed as before. The proof is similar to that of~\eqref{eq:thm_online_1}, so we only show the first stage that compares $\pi_\vartriangle$ with $\pi_l$, the optimal policy for $Q_l = (c_l,r,W)$ with $c_l(z) = K_l \1(z) + lz$: 
    \begin{equation*}
        \frac{\frac{1}{3M} J_{\pi_\vartriangle}(x_0 \given Q_h)}{J_{\pi_l}(x_0 \given  Q_l)}
        \hintedrel[aff_onl]{\leq} \frac{J_{\pi_l}(x_0 \given Q_h)}{J_{\pi_l}(x_0 \given  Q_l)}
        \leq \max \left\{ \frac{K_h}{K_l},\frac{h}{l} \right\},
    \end{equation*}
    where $Q_h = (c_h,r,W)$ with $c_h(z) = K_h \1(z) + hz$, and we use~\eqref{eq:online_affine} for (\hintref{aff_onl}). All of the inequalities hold for any $x_0$ and $Q$ with the appropriate ordering cost function. Therefore, they hold for the supremum over $x_0$ and $Q$. 
\end{proof}
Proposition~\ref{prop:online_algs} establishes that one can apply the online policy $\pi_\vartriangle$ to multi-location inventory problems with a variety of ordering cost functions and achieve bounded suboptimality. While these policies may be far from optimal in general, we show in the following section that they perform reasonably well for representative numerical examples, suggesting that one need not solve a dynamic program to achieve good performance in multi-location settings with nonlinear costs.

\section{Simulations}\label{sec:sims}
In this section we evaluate the suboptimality of the proposed decoupled policies and online algorithms via numerical simulation. We use the sector-bounded ordering cost function
\begin{equation}\label{eq:sim_sector}
    c(z) = \begin{cases}
        4z, \quad & z \leq 6 \\
        2z + 12, \quad & z > 6
    \end{cases}
\end{equation}
and the affine-bounded ordering cost function
\begin{equation}\label{eq:sim_affine}
    c(z) = \begin{cases}
        4\1(z) + 2z, \quad & z \leq 6 \\
        4\1(z) + z + 6, \quad & z > 6
    \end{cases}
\end{equation}
Both these functions exhibit economies of scale where the per-unit cost decreases as larger quantities are ordered. We perform $1000$ Monte Carlo simulations from each possible initial condition to compare decoupled dynamic programming-based policies and online algorithms with the optimal policy. 

\begin{figure*}
    \centering
    \includegraphics[width=\textwidth]{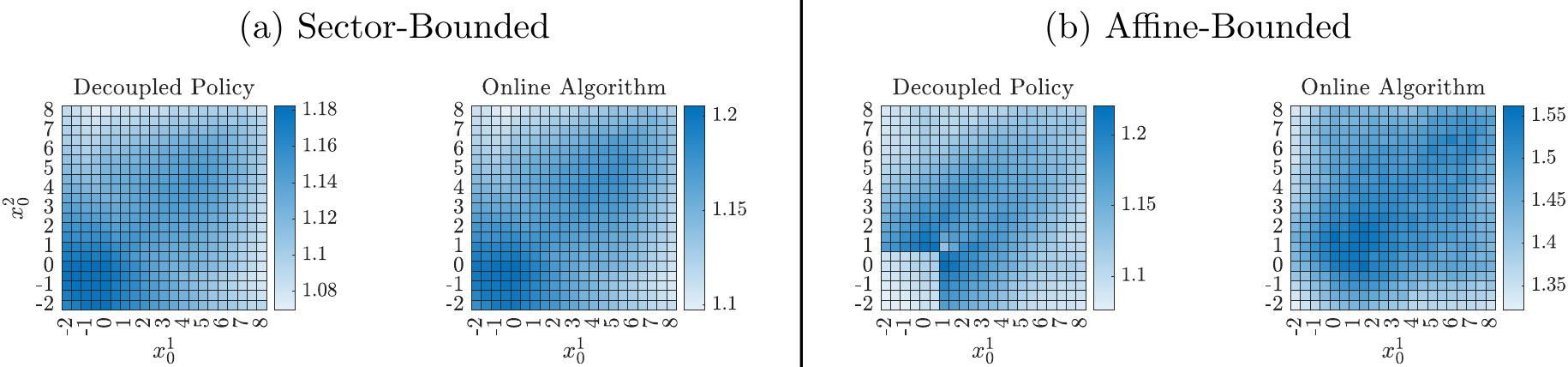}
    \caption{For all simulations we consider $M = 2$ inventories, state space $x \in [-2,8]$, discretization $\Delta x = 0.5$, a finite horizon $N = 20$, and stochastic demand $w_k \in \{0,0.5,1,1.5\}$ where $w_k \sim \mathrm{Bin}(3,0.5)$. (a) shows average cost ratios for the decoupled policy $\pi_\Hsquare$ and online policy $\pi_\vartriangle$ from each initial state under sector-bounded ordering cost. We set holding/backlog cost $r(z) = 0.1 \max\{0,z\} + 10 \max \{0,-z\}$. (b) shows average cost ratios for the decoupled policy $\pi_s$ and online policy $\pi_\vartriangle$ from each initial state under affine-bounded ordering cost. We set holding/backlog cost $r(z) = 0.2 \max\{0,z\} + 10 \max \{0,-z\}$.}
    \label{fig:sim}
\end{figure*}

Fig.~\ref{fig:sim}(a) shows a heatmap of the average cost ratios $J_{\pi_\Hsquare}(x_0 \given Q)/J_{\pi_*}(x_0 \given Q)$ and $J_{\pi_\vartriangle}(x_0 \given Q)/J_{\pi_*}(x_0 \given Q)$ under the ordering cost~\eqref{eq:sim_sector} where $\pi_\Hsquare$ is the decoupled base-stock policy that is optimal for the problem $P_l = (c_l,r,W)$ with $c_l(z) = lz$. For the decoupled policy, the mean ratio is 1.13 and the max ratio is 1.18. For the online algorithm, the mean ratio is 1.15 and the max ratio is 1.21. Our worst-case bounds are $h/l = 4/2 = 2$ for the decoupled policy and $2h/l = 4$ for the online algorithm. While the empirical performance is far from the worst-case bounds, we do see substantially increased cost. Interestingly, for this choice of parameters, the online algorithm performs nearly as well as the decoupled policy. 

Fig.~\ref{fig:sim}(b) shows a heatmap of the average cost ratios $J_{\pi_\diamond}(x_0 \given Q)/J_{\pi_*}(x_0 \given Q)$, $J_{\pi_\vartriangle}(x_0 \given Q)/J_{\pi_*}(x_0 \given Q)$ under~\eqref{eq:sim_affine} where $\pi_\diamond$ is the decoupled $(s,S)$ policy that is optimal for each single-location problem $P_h^i = (c_h,r^i,W^i)$ with $c_h(z) = K_h \1(z) + hz$. For the decoupled policy, the mean ratio is 1.14 and the max ratio is 1.22. For the online algorithm, the mean ratio is 1.47 and the max ratio is 1.56. Our (optimized) worst-case bounds are $M\max \{K_h/K_l,h/l\} = 2\max \{(16/3)/4,(4/3)/1\} = 8/3$ and $3M \max\{K_h/K_l,h/l\} = 8$, respectively. Once again the empirical performance does not attain the worst-case bounds; however, the suboptimality of both policies is evident. 
    
\section{Conclusion}
In this work, we evaluated the performance of simple decoupled policies in multi-location inventory control problems coupled through a nonlinear ordering cost function. For sector-bounded nonlinearities, we demonstrated that a collection of base-stock policies could achieve a tight bound parameterized by the nonlinearity of the function, while for affine-bounded nonlinearities, we demonstrated that a collection of $(s, S)$ policies could achieve a bound parameterized by the nonlinearity of the function and the number of inventories. 
%\james{tightness reminder}
We also showed that online algorithms from the literature could be used as a simpler alternative to dynamic programming in these settings. Numerical simulations indicate that for many problem instances there is additional structure that leads to better performance than our worst-case bounds
%\james{Preceding statement has a more negative tone than necessary. It's not surprising that the simulations perform much better than the bounds, this is standard for worst-case analyses. The positive conclusion is that practical problem instances have some additional `niceness' structure that our sector/affine/etc. classes fail to capture. Figuring out the right structure should give us better instance-dependent upper bounds that are closer to the simulations. You touch on this in the following sentence, but it could be motivated more strongly.}
. Future work includes deriving corresponding lower bounds for the performance of the decoupled policies as well as identifying particular problem classes that could admit tighter bounds.

\bibliographystyle{plain}
\bibliography{bib/New}

\appendix
\subsection{Cost Transformation}\label{app:cost_transformation}
Here we recap a useful cost transformation from the inventory control literature~\cite{levi2007approximation,janakiraman2004costs}. Consider rewriting the ordering cost $c(u_k)$ as $c(x_{k+1} - x_k + w_k)$ using the dynamics~\eqref{eq:dynamics}. For stationary ordering cost functions $c$ with a linear term $m \sum_{i=1}^M u_k^i$, the transformed cost $\hat{c}$ that neglects the linear term, i.e., sets $m = 0$, can be seen to neglect a cost of $m \sum_{k=0}^{N-1} \sum_{i=1}^M  w_k^i$ which is only a function of the stochastic demand. Therefore, we have the following relationship between the expected cost of a given policy $\pi$ under the original costs $P = (c,r,W)$ and the transformed costs $\hat{P} = (\hat{c},r,W)$:
\begin{equation}\label{eq:cost_transformation}
    J_\pi(x_0 \given P) {=} J_\pi(x_0 \given \hat{P}) + m \E \left[ \frac{1}{N} \sum\nolimits_{k=0}^{N-1} \sum\nolimits_{i=1}^M w_k^i \right].
\end{equation}

\subsection{Tightness Example for Theorem \ref{thm:main}}
\begin{proof}\label{app:tightness_thm_1}
    Consider an $M$-location infinite horizon problem $P$ with demand $w_k^i \sim \mathrm{Unif}(1, 1 + \delta)$ for all $i, k$ and holding/backlog cost $r^i(z_k^i) = \delta \max \{0, z_k^i \} + p \max \{0, -z_k^i \}$ where $\delta = \frac{\varepsilon}{l + 2}$. We define the set $V = \{M, M(1 + \delta) \}$. The ordering cost is $c(z) = l z$ for $z \in V$, and $c(z) = hz$ otherwise; we assume $p \gg h$. 
    
    We begin by deriving a lower bound on the expected cost of a decoupled base-stock policy $\pi_b \in \Pi_\textup{B}^M$ from the initial condition where $x_0^i = 0$ for all $i$. The policy is arbitrary aside from the assumption that $\sum_{k=0}^{n} S_k^i \geq n + 1$ (equivalently, $\sum_{k=0}^{n} u_k^i \geq n + 1$) for all $i$ and for all $n \geq 0$; it is easily verifiable that any policy that does not satisfy this condition is suboptimal since it incurs backlog costs with probability $1$. We can rewrite the dynamics as
    \begin{equation*}
        x_{k + 1}^i = \begin{cases}
            S_k^i - w_k^i & x_k^i < S_k^i \\
            x_k^i - w_k^i & x_k^i \geq S_k^i,
        \end{cases}
    \end{equation*}
    from which it is clear that the state $x_k^i$ is distributed as the sum of at most $k$ $\mathrm{Unif(0, \delta)}$ random variables (with mean possibly shifted by the base stock levels $S_0, ..., S_k$) for all $k \geq 1$. Since $V$ is a measure-zero subset of $\mathbb{R}$, it follows that
    \begin{equation*}
        \mathbb{P} \left[ \sum\nolimits_{i = 1}^M \max \{S_k^i - x_k^i, 0 \} \in V \right] = 0,
    \end{equation*}
    so every unit of inventory ordered according to $\pi_b$ for all $k \geq 1$ incurs an ordering cost of $h$ with probability $1$. Thus, for any finite $N \geq 1$,
    \begin{align*}
        J_{\pi_b}(x_0 \given P ) 
        &\geq \mathbb{E} \left[ \sum\nolimits_{k = 0}^{N - 1}  c(u_k) \right]
        % = \mathbb{E} \left[ c(u_0) + \sum\nolimits_{k = 1}^{N - 1}  c(u_k) \right]
        \\
        &\geq \sum\nolimits_{i = 1}^M l \cdot u_0^i +\sum\nolimits_{k = 1}^{N - 1} \sum\nolimits_{i = 1}^M h \cdot u_k^i \\
        &\geq Ml + (N-1)Mh.
    \end{align*}
    
    Next, we upper bound the cost of a policy we call $\pi_v$ which simply places an order for the minimum value of $v \in V$ such that $v + \sum_{i = 1}^M x_k^i \geq M(1 + \delta)$ and allocates it so that the inventory level in each location is equal. It is clear that $\pi_v$ incurs an ordering cost of $l$ per unit of inventory, and it holds no more than $2 \delta$ units of inventory in each location in each timestep. For finite $N$, this yields
    \begin{align*}
        J_{\pi_v}(x_0 \given P ) &= \E \left[ \sum\nolimits_{k = 0}^{N - 1} c(u_k) + r(x_k + u_k - w_k) \right] \\
        &\leq \sum\nolimits_{k = 0}^{N - 1} l \cdot M (1 + \delta) + \sum\nolimits_{k = 0}^{N - 1} \sum\nolimits_{i = 1}^{M} 2 \delta \\
        &\leq N M (1 + \delta) l + 2 N M \delta.
    \end{align*}
    Since $\pi_v$ is not necessarily $\pi_*$, we have that as $N \to \infty$,
    \begin{align*}
        \limsup_{N \to \infty} \frac{J_{\pi_b}(x_0 \given P )}{J_{\pi_*}(x_0 \given P )}
        &\geq \limsup_{N \to \infty} \frac{J_{\pi_b}(x_0 \given P )}{J_{\pi_v}(x_0 \given P )} \\
        &\geq \limsup_{N \to \infty}  \frac{MNh + M(l - h)}{M N (l + \delta(l + 2))} \\
        &= \frac{h}{l + \delta(l + 2)} = \frac{h}{l + \varepsilon}
    \end{align*}
    where the last line follows from the fact that the term $M(l - h)$ is a fixed constant. One can construct limiting problem instances by taking $\varepsilon \to 0$ to yield the supremum $\frac{h}{l}$.

    % To extend this result to the infinite horizon and to any initial condition, we modify $\pi_v$ to order $\sum_{i = 1}^M u_k^i = \max\{1 + \delta - x_k^i, 0\}$ units of inventory if any $x_k^i \in (-\infty, -\delta] \ \cup \ [\delta, 1 + \delta)$. Thus, in exactly one timestep, $x_k^i \geq 1 + \delta$ for all $i$. For any location where $x_0^i \geq 1 + \delta$, it follows that for some $k \leq k^i = \lfloor x_0^i \rfloor$,  $x_k^i \leq \delta$. Hence, for all timesteps $k \geq k' = \max_i k^i$, the behavior of the system under $\pi_v$ from any initial condition resembles its behavior for initial conditions $(-\delta, 0)$ as described above. Furthermore, for all $k \geq 0$, the lower bound on the cost of $\pi_b$ remains the same. The cumulative costs of both $\pi_b$ and $\pi_v$ for $k \leq k'$ are finite, so the ratio of the average costs as $N \to \infty$ is the same as in the analysis above, yielding $\frac{J_{\pi_b}(x_0 \given P)}{J_{\pi_v}(x_0 \given P)} \geq \frac{h}{l + \varepsilon}$.
\end{proof}

\end{document}